\RequirePackage{fix-cm}
\documentclass[smallcondensed]{svjour3}       
\smartqed  
\usepackage{graphicx}
\usepackage{mathptmx}      
\usepackage{amscd, amssymb, amsfonts, amsmath}
\usepackage[all]{xy}
\usepackage{color}
\usepackage{mathdots}
\usepackage{hyperref}
\usepackage{latexsym}
\newcommand{\gl}{\mbox{\rm{gldim}}}
\newtheorem{teo}{Theorem}[section]
\newtheorem{prop}[teo]{Proposition}
\newtheorem{lema}[teo] {Lemma}
\newtheorem{ej}[teo]{Example}
\newtheorem{obs}[teo]{Remark}
\newtheorem{defi}[teo]{Definition}

\DeclareMathOperator{\len}{l}
\DeclareMathOperator{\st}{s}
\DeclareMathOperator{\tg}{t}
\DeclareMathOperator{\fidim}{\phi dim}

\DeclareMathOperator{\md}{mod}
\DeclareMathOperator{\ind}{ind}
\DeclareMathOperator{\add}{add}

\DeclareMathOperator{\pd}{pd}
\DeclareMathOperator{\id}{id}

\DeclareMathOperator{\fin}{findim}
\DeclareMathOperator{\rd}{repdim}

\DeclareMathOperator{\ext}{Ext}
\DeclareMathOperator{\Ker}{Ker}
\DeclareMathOperator{\Tp}{Top}

\begin{document}

\title{On algebras of $\Omega^n$-finite and $\Omega^{\infty}$-infinite representation type}


\author{Marcos Barrios \and Gustavo Mata}


\institute{M. Barrios \at
              Universidad de la Rep\'ublica, Montevideo, Uruguay\\
              \email{marcosb@fing.edu.uy}           
           G. Mata \at              Universidad de la Rep\'ublica, Montevideo, Uruguay\\
              \email{gmata@fing.edu.uy}}

\maketitle

\begin{abstract}

Co-Gorenstein algebras were introduced by A. Beligiannis in \cite{B}. In \cite{KM}, the authors propose the following conjecture (Co-GC): if  $\Omega^n (\md A)$ is extension closed for all $n \leq 1$, then $A$ is right Co-Gorenstein, and they prove that the Generalized Nakayama Conjecture implies the Co-GC, also that the Co-GC implies the Nakayama Conjecture. In this article we characterize the subcategory $\Omega^{\infty}(\md A)$ for algebras of $\Omega^{n}$-finite representation type. As a consequence, we characterize when a truncated path algebra is a Co-Gorenstein algebra in terms of its associated quiver. We also study the behaviour of Artin algebras of $\Omega^{\infty}$-infinite representation type. Finally, it is presented an example of a non Gorenstein algebra of $\Omega^{\infty}$-infinite representation type and an example of a finite dimensional algebra with infinite $\phi$-dimension.

\end{abstract} 

\keywords{Truncated path algebra, Co-Gorenstein algebra, Igusa-Todorov function.}
\PACS{16W50\and 16E30}
\subclass{16G10}

\section{Introduction}

The Finitistic Dimension Conjecture (FDC\footnote{$\fin(A) = \sup\{\pd (M)|M \in \md A
\text{ with } \pd(M) < \infty\} < \infty.$}), the Generalized Nakayama Conjecture (GNC\footnote{For any simple module $S \in \md A$, there exists $i \geq 0$ such that $\ext^i_A (S, A)  = 0.$}) and the Nakayama Conjecture (NC\footnote{Let $0 \rightarrow A \rightarrow I^0 \rightarrow I^1 \rightarrow \ldots$ be a minimal injective resolution of the $A$-module $A$. If $I^i$ is projective for any $i \geq 0$, then $A$ is self-injective.}) are examples of important conjectures in Representation Theory of Artin algebras and they still remain open. These conjectures are related in the following way:

If $A$ is an Artin algebra,
\begin{itemize}

\item If $A$ verifies the FDC then $A$ verifies the GNC,

\item If $A$ verifies the GNC then $A$ verifies the NC.

\end{itemize}

An algebra $A$ is Co-Gorenstein if the subcategories $\Omega^{\infty}(\md A)$(Definition \ref{sizigia infinita}) and $\mathcal{G}\mathcal{P}(A)$ (Section 2) agree. In \cite{KM} the authors prove that there is a relation between NC, GNC and Co-Gorenstein algebras. They prove that GNC implies Conjecture 1 (if $\Omega^n(\md A)$ is extension closed for all $n \geq 1$, then $A$ is right Co-Gorenstein), and they also prove that Conjecure 1 implies NC.
This justifies the relevance of studying this family of algebras. 

The main result of this article is the following.\\

\underline{{\bf Theorem 3.9}}
For an Artin algebra $A$ of $\Omega^n$ finite representation type, a module is in the subcategory
$\Omega^\infty(\md A)$ if and only if it is either projective or $\Omega$-periodic.\\


In particular, truncated path algebras are of $\Omega^1$-finite representation type (Theorem 5.11 of \cite{BZR}). We use both theorems and determine, in terms of the geometry of the quiver, when a truncated path algebra is a Co-Gorenstein algebra (Theorem \ref{clasificacion coGorenstein}).

In an attempt to prove the FDC, K. Igusa and G. Todorov define, in \cite{IT}, two functions from the objects of $\md A$ to the natural numbers, which generalizes the notion of projective dimension. Nowadays, they are known as the Igusa-Todorov functions, $\phi$ and $\psi$. In that article the authors prove, with the Igusa-Todorov functions, that if an Artin algebra $A$ has $\rd(A) \leq 3$ then $A$ verifies the FDC. 

More recently, it was conjectured in \cite{FLM} that the associated dimensions to these functions are finite ($\phi$DC\footnote{$\fidim(A) = \sup \{ \phi(M): M\in \md(A) \}< \infty$}). This conjecture is also related with the FDC as follows: For $A$ an Artin algebra, if $A$ verifies the $\phi$DC then $A$ verifies the FDC. In \cite{LM}, the authors prove that every Artin algebra of $\Omega^n$-finite representation type and every Gorenstein algebra verifies $\phi$DC. Therefore, to find an example of an algebra that does not verify the $\phi$DC, it should be looked in the family of non Gorenstein algebras of $\Omega^{\infty}$-infinite representation type. 

In this article we show that there are algebras of $\Omega^{\infty}$-infinite representation type that are not Gorenestein algebras (Examples \ref{finito} and \ref{infinito}). Finally we exhibit an example (Example \ref{infinito}) where $\phi$DC is not true, which implies that $\phi$DC and FDC are not equivalent. We note
that $\phi$DC is disproven independently by both the present paper and \cite{HI}, with
preliminary versions of both papers appearing on the arXiv in November 2019.

\section{Preliminaries}

\subsection{Notations and terminology}

From now on $A$ will always be an Artin algebra. We denote by $\md A$ the category of finitely generated right $A$-modules. The full subcategory of $\md A$ whose objects are direct summands of finite sums of copies of $M$ is denoted by $\add (M)$. For a module $M \in \md A$, $\Tp (M)$, $\Omega (M)$ and $\pd (M)$ are the top, the syzygy and the projective dimension of $M$, respectively. For an algebra $A$, $\gl (A)$ and $\fin (A)$ are the global and finitistic dimension of $A$, respectively. 

We denote by $\mathcal{P}(A)$ the full subcategory formed by projective modules and by $\mathcal{G}\mathcal{P}(A)$ the full subcategory formed by Gorenstein projective modules. Observe that $\mathcal{P}(A) \subset \mathcal{G}\mathcal{P}(A)$.
Recall that an algebra $A$ is self-injective if and only if $\mathcal{G}\mathcal{P}(A) = \md A$, that is, all modules are Gorenstein projective. On the other hand, $A$ is called CM-free provided that $\mathcal{G}\mathcal{P}(A) = \mathcal{P}(A)$, that is, all its finitely generated Gorenstein projective modules are projective. Recall that for an Artin algebra $A$, the set $\ind \  \underline{\mathcal{G}\mathcal{P}(A)}$ of isoclasses of indecomposable objects in the stable category $\mathcal{G}\mathcal{P}(A)$ is identified with the set of isoclasses of indecomposable non-projective Gorenstein-projective $A$-modules.

We denote by $\mathcal{P}^{< \infty} (A)$ the full subcategory of $\md A$ formed by the modules of finite projective dimension and by $\mathcal{P}^{\infty} (A)$ the full subcategory of $\md A$ formed by the modules of infinite projective dimension with no summand with finite projective dimension.

Given $n \in \mathbb{N}$, we denote by $\add (\Omega^n (\md A))$ the full subcategory of $\md A$ formed by the direct sums of direct summands of modules in $\Omega^n (\md A)$. For $n \in \mathbb{N} $, we say that $\Omega^n (\md A)$ is of {\bf{finite representation type}} if there are finitely many indecomposable modules $M$ for which there exists $N \in \md A$ with $M$ a direct summand of $\Omega^n (N)$. In this case $A$ is called of {\bf $\Omega^n$-finite representation type} (for short {\bf syzygy finite}). If $A$ is not of $\Omega^n$-finite representation
type for any $n \in \mathbb{N}$, we say that $A$ is {\bf $\Omega^{\infty}$-infinite representation type}. 

Notice that if $A$ is a syzygy finite algebra then $\fin(A) < \infty$.  

\subsection{Quivers and path algebras}

If $Q = (Q_0,Q_1,\st,\tg)$ is a finite connected quiver, $\Bbbk Q$ denotes its associated path algebra. Given $\rho$ a path in $\Bbbk Q$, we denote by $\len(\rho)$, $\st(\rho)$ and $\tg(\rho )$ the length, start and target of $\rho$, respectively. $C^n$ denotes the oriented cycle graph with $n$ vertices. We denote by $J$ the ideal generated by all the arrows from $Q_1$. Given a quiver $Q$ with cycles, we denote by $Q^{\infty}$ the full subquiver of $Q$ formed by the vertices $w$ with infinite paths starting at $w$.
A full subquiver $Q'$ of a quiver $Q$ is called a \textbf{final subheart} (see \cite{LMM}) if it is minimal in the family of full subquivers of $Q$ which are closed by sucessors. We say a final subheart is {\bf trivial} if $Q'$ is a simple vertex. The proof of the following result is clear. 

\begin{obs} 
If $Q$ has a cycle and has no surces, then $Q$ has a non trivial final subheart.
\end{obs}

A finite dimensional algebra $A = \frac{\Bbbk Q}{I}$ is called a \textbf{monomial algebra} if $I$ is generated by paths.

Let $A = \frac{\Bbbk Q}{I}$ be an algebra, if $v \in Q_0$ we denote by $S_v$ and $I_v$ the simple $A$-module and the injective $A$-module associated to $v$, respectively. For $M \in \md \frac{\Bbbk Q}{I}$, we denote by $\underline{\dim}(M)$ its dimension vector.

Let $A = \frac{\Bbbk Q}{I}$ be a monomial algebra. We call a pair $(p, q)$ of nonzero paths in A a \textbf{perfect pair} provided that the following conditions are
satisfied:

\begin{itemize}

\item both of the nonzero paths $p$, $q$ are nontrivial satisfying $\st(p) = \tg(q)$ and
$qp = 0$ in $A$;

\item if $q'p = 0$ for a nonzero path $q'$ with $\tg(q') = \st(p)$, then $q'= q''q$ for some
path $q''$;

\item if $qp'= 0$ for a nonzero path $p'$ with $\st(p') = \tg(q)$, then $p' = pp''$ form some
path $p''$;

\end{itemize}

Let $A = \frac{kQ}{I}$ be a monomial algebra. We call a nonzero path $p$
in $A$ a \textbf{perfect path}, provided that there exists a sequence:
$$p = p_1, p_2 , \ldots , p_n , p_{n+1} = p$$
of nonzero paths such that $(p_{i} , p_{i+1})$ are perfect pairs for all $1 \leq i \leq n$. If
the given nonzero paths $p_i$ are pairwise distinct with the exception that $p_1 = p_{n+1}$, we refer to the sequence $p =
p_1 , p_2 , \ldots , p_n , p_{n+1} = p$ as a \textbf{relation-cycle} for $p$.

The next theorem shows how to compute every non projective Gorenstein projective module. 

\begin{teo}\cite[Theorem 4.1]{CSZ}\label{CSZ}
Let $A$ be a monomial algebra. Then there is a bijection
$$\{\mbox{perfect paths in } A\} \rightarrow \ind\ \underline{\mathcal{G}\mathcal{P}(A)}$$
sending a perfect path $p$ to the A-module $Ap$.
\end{teo}

The following theorem shows that monomial algebras are of $\Omega^2$-finite representation type algebras.

\begin{teo}\cite[Theorem 3.I]{Z}\label{Z}
Let $A$ be a finite dimensional monomial algebra. If $M$ is an $A$-module such that $M \cong \Omega^t(N)$ for an $A$-module $N$ and $t \geq 2$, then it is a direct sum of right ideals of $A$ which are generated by paths with length greater or equal than $1$. In particular, monomial algebras are of $\Omega^2$-finite representation type.
\end{teo}

\subsection{Truncated path algebras}

A monomial algebra $A = \frac{\Bbbk Q}{I}$ is called a \textbf{truncated path algebra} if there is an $m \in \mathbb{N}, m\geq 2$ such that $I = J^m$. If $m = 2$ we say that $A$ is a \textbf{radical square zero algebra}. Following \cite{BMR}, we use the below notation.

For a truncated path algebra $A$, we denote by $ _\rho M^l_v(A)$ the ideal $\rho A$, where $\len(\rho) = l$, $\tg(\rho) = v$

\begin{obs}
Notice that if $\rho$ and $\mu$ are paths such that $\tg(\rho) = \tg(\mu) = v$ and $\len(\rho) = \len(\mu) = l$ then $ _\rho M^l_v(A) \cong  \hspace{0cm} _{\mu} M^l_v(A)$
\end{obs}

By the previous remark we can omit the path subindex so we write $M^l_v(A)$ instead of $ _\rho M^l_v(A)$. Finally we denote by $M^l(A) = \oplus_{v \in Q_0} M^l_v(A) $.

Given $0<l\leq m-1$, it is easy to see that $M^l_v(A)$ exists if and only if $\id_{\frac{\Bbbk Q}{J^{2}}} (S_v) \geq l$. Also, $M^l_v(A)$ is a projective $A$-module if and only if $\pd_{\frac{\Bbbk Q}{J^{2}}} (S_v) \leq m-l-1$. 

\begin{obs}\label{syzygy truncadas}
If $A = \frac{\Bbbk Q}{J^m}$ is a truncated path algebra, then
$$\Omega( M_v^l(A)) = \bigoplus_{\rho:\left\{\substack{ \st(\rho) = v\\ \len(\rho)= m-l}\right.} M_{ \tg(\rho)}^{m-l}(A),$$
\end{obs}

Proposition $2.6$ gives a method to compute the global dimension for truncated path algebras as a function of the length of the longest path of their quivers.

\begin{prop}\label{gldim truncadas}\cite[Theorem 5]{DH-ZL}
Let $A  \cong \frac{\Bbbk Q}{J^m}$ be a truncated path algebra. Then the following statements hold.
\begin{enumerate}

\item $\gl (A) < \infty$ if and only if $Q$ has no oriented cycles. 

\item If $A$ has finite global dimension, then 
$$\gl (A) =
\left\{
\begin{array}{ll}
2\frac{l}{m}&\text{if } l\equiv0\pmod{m},\\
2\left\lfloor\frac{l}{m}\right\rfloor+1&\text{otherwise},
\end{array}\right.$$ where $l$ is the length of the largest path of $\Bbbk Q$.

\end{enumerate}

\end{prop}

\subsection{Igusa-Todorov functions}

We recall the definition of the Igusa-Todorov $\phi$ function and some basic properties. We also define the $\phi$-dimension of an algebra.

\begin{defi}
Let $K_0(A)$ be the abelian group generated by all symbols $[M]$, where $M \in \md A$, modulo the relations
\begin{enumerate}
  \item $[M]-[M']-[M'']$ if  $M \cong M' \oplus M''$,
  \item $[P]$ for each projective.
\end{enumerate}
\end{defi}

\noindent Let $\bar{\Omega}: K_0 (A) \rightarrow K_0 (A)$ be the group endomorphism induced by $\Omega$
If $\mathcal{C}$ is subcategory of $\md A$ closed by direct sums and direct summands such that $\Omega(\mathcal{C}) \subset \mathcal{C}$, then we define by $K_0(\mathcal{C}) = \{[M]: M \in \mathcal{C}\}$. Now, if $M$ is a finitely generated $A$-module then $\langle \add M\rangle$ denotes the subgroup of $K_0 (A)$ generated by the classes of indecomposable summands of $M$.

\begin{defi}\label{monomorfismo}
\cite{IT} The \textbf{(right) Igusa-Todorov function} $\phi$ of $M\in \md A$  is defined as 
\[\phi_{A}(M) = \min\left\{l: \bar{\Omega}_{A} {\vert}_{{\bar{\Omega}_{A}}^{l+s}\langle \add M\rangle}\text{ is a monomorphism for all }s\in \mathbb{N}\right\}.\]
\end{defi}

In case that there is no possible misinterpretation we denote by $\phi$ the Igusa-Todorov function $\phi_A$.

\begin{prop}[\cite{IT}, \cite{HLM}] \label{it1} \label{Huard1}
Given $M,N\in$ mod$A$. 

\begin{enumerate}
  \item $\phi(M) = \pd (M)$ if $\pd (M) < \infty$.
  \item $\phi(M) = 0$ if $M \in \ind A$ and $\pd(M) = \infty$.
  \item $\phi(M) \leq \phi(M \oplus N)$.
  \item $\phi\left(M^{k}\right) = \phi(M)$ for $k \in \mathbb{N}$.
  \item $\phi(M) \leq \phi(\Omega(M))+1$.
\end{enumerate}
\begin{proof}
For the statements 1-4 see \cite{IT}, and for 5 see \cite{HLM}.
\end{proof}
\end{prop}

\begin{defi} Let $A$ an Artin algebra, and $\mathcal{C}$ a subcategory of $\md A$. The $\phi$-dimension of $\mathcal{C}$ is defined as
$$\fidim  (\mathcal{C}) = \sup \{ \phi(M): M \in \mathcal{C} \}.$$
In particular, we denote by $\fidim (A)$ the $\phi$-dimension of $\md A$
\end{defi}

The following theorem shows that the finitude of the $\phi$-dimension is invariant by derived equivalences.

\begin{teo} \cite[Theorem 4.10]{FLM}\label{FLM}
Let A and B be Artin algebras, which are derived equivalent. Then,
$\fidim (A) < \infty $ if and only if $\fidim (B) < \infty $.
\end{teo}

\section{Periodic modules over algebras of $\Omega^n$-finite representation type}

In this section we prove that the subcategory of $\Omega$-periodic $A$-modules and the subcategory $\Omega^{\infty}(\md A)$ agree for an algebra $A$ of $\Omega^n$-finite representation type.


\begin{defi}
Let $A$ be an Artin algebra. We define $\Vert \cdot \Vert : \md A \rightarrow \mathbb{N}$ a function such that $\Vert M \Vert = k$ if $M \cong (\oplus_{i=1}^k M_i) \oplus (\oplus_{j=1}^m N_j)$ is the descomposition into indecomposable modules where $\pd (M_i) = \infty$ $\forall i \in \{1, \ldots, k\}$ and $\pd (N_j) < \infty$ $\forall j \in \{1, \ldots, m\}$.
\end{defi}

\begin{obs} If $M$ and $N$ are $A$-modules, we have that:

\begin{enumerate}

\item $\Vert M \oplus N \Vert = \Vert M \Vert + \Vert N \Vert$.

\item $\Vert M^k \Vert = k \Vert M \Vert$ for $k \in \mathbb{N}$.

\item $\Vert M \Vert \leq \Vert \Omega (M)\Vert$. 

\item If $\Omega^m(M) \cong M$ then $\Vert \Omega (M)\Vert = \Vert M \Vert$.

\end{enumerate}

\end{obs}

\begin{defi}\label{sizigia infinita} For an algebra $A$ we denote by $\Omega^{\infty}(\md A)$ the full subcategory generated by the class of modules $M$ such that there exists an exact sequence  $0 \rightarrow M \rightarrow P_0 \rightarrow P_1 \rightarrow \ldots$  with $P_i \in \mathcal{P}(A)\  \forall  i \in \mathbb{N}$. We say that $\Omega^{\infty}(\md A)$ is {\bf trivial} if $\mathcal{P}(A) = \Omega^{\infty}(\md A)$. 
\end{defi}

\begin{defi}
When an $A$-module $M$ verifies $\Omega^m (M) \cong M$ for $m \in \mathbb{Z}^+$, we say that $M$ is {\bf $\Omega$-periodic}.
\end{defi}

\begin{obs} For an algebra $A$ we have the inclusions: 

$$\mathcal{P}(A) \subset \mathcal{G}\mathcal{P}(A) \subset \Omega^{\infty}(\md A).$$ 

If, in addition, the algebra $A$ has finite global dimension, then:

$$\mathcal{P}(A) = \mathcal{G}\mathcal{P}(A)$$

\end{obs}

\begin{obs} If the algebra $A$ has $\Omega^{\infty}(\md A)$ trivial, then $A$ is CM-free.

\end{obs}

\begin{obs}\label{infinitum}
If the algebra $A$ has $\fin (A) < \infty$, every non projective module $M$ in $\Omega^{\infty}(\md A)$ verifies $\pd (M) = \infty$.
\end{obs}

\begin{lema}\label{syzygy finite} If $A$ is of $\Omega^n$-finite representation type for $n \in \mathbb{Z}^{+}$, then the following assertions are equivalent

\begin{enumerate}

\item There is an $A$-module $M \neq \{0\}$ and $l \in \mathbb{Z}^{+}$ such that $\Omega^l (M) \cong M$.

\item There is a nonprojective module $\bar{M}$ in $\Omega^{\infty}(\md A)$.

\end{enumerate}

\begin{proof}

($1 \Rightarrow 2$) Suppose there exist $M \in \md A$ and $l \in \mathbb{Z}^+$ such that $\Omega^l (M) \cong M$. Then there exist an exact sequence
$$\zeta: \ 0 \rightarrow M \rightarrow P_l \rightarrow \ldots \rightarrow P_1 \rightarrow P_0 \rightarrow M \rightarrow 0.$$
We can join infinite exact sequences $\zeta$, thus we conclude that $M \in \Omega^{\infty}(\md A)$.

($2 \Rightarrow 1$) Suppose that $A$ is of $\Omega^n$-finite representation type and there exists a nonprojective module $\bar{M}$ in $\Omega^{\infty}(\md A)$. By Remark \ref{infinitum} $\bar{M} = M_1 \oplus M_2$ with $M_1 \in \mathcal{P}^{\infty}(A)$ and $M_2 \in \mathcal{P}^{< \infty}(A)$. Consider the set $\mathcal{W}_{\bar{M}} = \{ N \in \mathcal{P}^{\infty}(A) \cap \add ( \Omega^l(\md A) ) \mbox{ such that } \Vert N \Vert \leq \Vert M \Vert \}$. 
Clearly $\mathcal{W}_{\bar{M}}$ is a finite non empty set ($M_1 \in \mathcal{W}_{\bar{M}}$), thus there exist $N_0 \in \mathcal{W}_{\bar{M}}$, $f(k) = n_k : \mathbb{Z}^{+} \rightarrow \mathbb{Z}^{+}$ an increasing function, and  $N_{n_k} \in \mathcal{P}^{<\infty}(A)$ such that $\Omega^{n_k} (N_0 \oplus N_{n_k}) \cong {\bar{M}}$ for all $k \in \mathbb{N}$. 
We deduce that there are $i,j \in \mathbb{Z}^{+}$ such that $n_j - n_i > m = \fin (A)$ and $\Omega^{n_j-n_i}(N_0) \cong N_0 \oplus N'$ with $N' \in \mathcal{P}^{< \infty}(A)$. Finally we conclude that $N_0\oplus N'$ is a $\Omega$-periodic module.\end{proof}

\end{lema}

In the proof of the lemma above we observe that $\bar{M} \cong \Omega^k(N)$ for $N$ a periodic module and $k > \fin (A)$. So we deduce the next result. 

\begin{teo}\label{equivalencia periodicos} Let $A$ be an Artin algebra of $\Omega^n$-finite representation type. If $M$ is a non projective module, then $M \in\Omega^{\infty}(\md A)$ if and only if there exist $m \in \mathbb{Z}^+$ such that $\Omega^m(M) \cong M$.

\end{teo}

By Theorem \ref{Z}, Monomial algebras verify the hypoteses of Lemma \ref{syzygy finite} and Theorem \ref{equivalencia periodicos}.  

%
%
%
%



The following example shows that Theorem \ref{equivalencia periodicos} is not valid if the algebra is not syzygy finite.

\begin{ej} \label{ej sizigia-infinita} Let $\Bbbk$ be an infinite field with characteristic different to $2$. Consider the algebra $A = \frac{\Bbbk Q}{•I}$, where $Q$ is the quiver below:

$$ Q = 
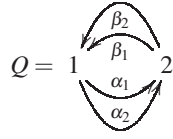
\xymatrix{ 1 \ar@/_{1pc}/[r]^{\alpha_1} \ar@/_{2pc}/[r]^{\alpha_2} & 2 \ar@/_{1pc}/[l]^{\beta_1} \ar@/_{2pc}/[l]^{\beta_2} }$$

and $I = \langle \beta_1\alpha_1-2\beta_2\alpha_2, \alpha_1\beta_1 - \alpha_2\beta_2,  \alpha_i \beta_j , \beta_i \alpha_j \mbox{ with } i,j \in \{1,2\}, i\not = j\rangle$.

The following assertions are clear.
\begin{itemize}

\item $ M_a = \xymatrix{ \Bbbk \ar@/_{1pc}/[r]^{a1_{\Bbbk}} \ar@/_{2pc}/[r]^{1_{\Bbbk}} & \Bbbk \ar@/_{1pc}/[l]^{0} \ar@/_{2pc}/[l]^{0} }$ are indecomposable modules for $a \neq 0$. 

\item $ N_a = \xymatrix{ \Bbbk \ar@/_{1pc}/[r]^{0} \ar@/_{2pc}/[r]^{0} & \Bbbk \ar@/_{1pc}/[l]^{1_{\Bbbk}} \ar@/_{2pc}/[l]^{a1_{\Bbbk}} }$ are indecomposable modules for $a \neq 0$.

\item $M_a \ncong M_b$ and $N_a \ncong N_b$ if $a \neq b$.

\item $\Omega (M_a) = N_{-a}$ and $\Omega (N_{a}) = M_{-\frac{a}{2}}$ for $a \neq 0$.

\end{itemize}

Then, we deduce that the modules $M_a \in \Omega^{\infty}(\md A)$ but are not $A$-periodic.  

\end{ej}

\section{Co-Gorenstein algebras}

The notion of Co-Gorenstein algebra was introduced by Beligiannis in \cite{B}. It is somewhat dual to the definition of Gorenstein algebra, whence the name (see Remark 5 of \cite{KM}). Now we characterize when a truncated path algebra is Co-Gorenstein.

\begin{defi} An algebra $A$ is right Co-Gorenstein if and only if one of the following equivalent conditions hold

\begin{itemize}

\item $\Omega^{\infty}(\md A) \subset  ^{\perp} \! \mathcal{P}(A)$.

\item $\Omega^{\infty}(\md A) = \mathcal{G}\mathcal{P}(A)$.

\end{itemize}

where $^\perp \mathcal{P}(A) =\{M\in \md A:\ \ext^i(M,A)=0, \ \forall i\geq 1\}$.

\end{defi}

The remark below follows from Corollary 3.4 in \cite{Zh}.

\begin{obs} Gorenstein algebras are Co-Gorenstein

\end{obs}

Using Remark \ref{syzygy truncadas}, we obtain a result that allows us compute the modules in $\Omega^{\infty}(\md A)$ for a truncated algebra $A$.

\begin{prop}\label{periodico} Let $A = \frac{\Bbbk Q}{J^m}$ be a truncated path algebra with infinite global dimension, then the following are equivalent

\begin{enumerate}

\item There exists an $A$-module $M \cong M_1 \oplus M_2 \in \Omega^{\infty}(\md A)$ where\\ $M_1 \cong \bigoplus_{i=1}^k M^{l_i}_{v_i}(A) \in \mathcal{P}^{\infty}(A)$ and $M_2 \in \mathcal{P}^{< \infty}(A).$
\item There exists a $\frac{\Bbbk Q^{\infty}}{J^m}$-module $N \cong N_1 \oplus P \in \Omega^{\infty}(\md \frac{\Bbbk Q^{\infty}}{J^m})$ where\\ $N_1 \cong \bigoplus_{i=1}^k M^{l_i}_{v_i}(\frac{\Bbbk Q^{\infty}}{J^m}) \in \mathcal{P}^{\infty}(\frac{\Bbbk Q^{\infty}}{J^m})$ and $P \in \mathcal{P}(\frac{\Bbbk Q^{\infty}}{J^m})$.
\end{enumerate}

\end{prop}

The lemma below follows from Remark \ref{syzygy truncadas}.

\begin{lema}\label{norma que crece}
Let $Q$ be a quiver such that $Q^{\infty} \not= \emptyset$. Let $v$ be a vertex in $Q^{\infty}$ such that there are two paths in $Q^{\infty}$ $\rho$ and $\rho'$ starting at $v$ where $\rho \not \subset \rho'$ and $\rho' \not \subset \rho$. Then there exist $m_v^l \in \mathbb{Z}^+$ such that $\Vert M_v^l(A) \Vert < \Vert \Omega^{m_v^l}( M_v^l(A))\Vert $.
\end{lema}

The following result characterize when a truncated path algebra is a Co-Gorenstein algebra.

\begin{teo}\label{clasificacion coGorenstein} If $A = \frac{\Bbbk Q}{J^m}$ is a truncated path algebra, then the following are equivalent:

\begin{enumerate}

\item $A$ is Co-Gorenstein.

\item $Q$ is one of the following types:

\begin{itemize}

\item[i] $Q$ is acyclic.

\item[ii] $Q = C^n$.

\item[iii] $Q^{\infty}$ has no final subheart equal to $C^i$ for $i \in \mathbb{Z}^+$.

\end{itemize}  

\end{enumerate}

\begin{proof}

Any finite quiver $Q$ belongs to only one of the following families: 

\begin{enumerate}

\item[i] $Q$ acyclic,

\item[ii] $Q = C^n$,

\item[iii] $Q$ has cycles, $Q \neq C^n$ and every final subheart of $Q^{\infty}$ is not equal to $C^i$ for $i \in \mathbb{Z}^+$,

\item[iv] $Q$ has cycles, $Q \neq C^n$ and there exist a final subheart of $Q^{\infty}$ equal to $C^i$ for $i \in \mathbb{Z}^+$.

\end{enumerate}

\underline{Claim:} if $Q$ is in the first three families, then $A$ is Co-Gorenstein.

\begin{itemize}

\item[i] If $Q$ is acyclic, then $\gl A < \infty$. Thus $A$ is a Gorenstein algebra.

\item[ii] If $Q = C^n$, then $A$ is a self-injective algebra. Thus $A$ is $0$-Gorenstein. 

\item[iii] By Lemma \ref{norma que crece}, for every module $M_v^l(\frac{\Bbbk Q^{\infty}}{J^m})$ there exist $m_v^l \in \mathbb{Z}^+$ such that $\Vert M_v^l(A) \Vert < \Vert \Omega^{m_v^l}( M_v^l(A))\Vert $. Using the last fact and that every module in $ \Omega^{\infty}(\frac{\Bbbk Q^{\infty}}{J^m})$ must be a direct sum of modules $M_v^l(\frac{\Bbbk Q^{\infty}}{J^m})$ and a projective module we obtain the equality $ \Omega^{\infty}(\frac{\Bbbk Q^{\infty}}{J^m}) = \mathcal{P}(\frac{\Bbbk Q^{\infty}}{J^m})$.
By Corollary \ref{equivalencia periodicos} we have that there are no $\Omega$-periodic $\frac{\Bbbk Q^{\infty}}{J^m}$-modules. By Proposition \ref{periodico} we have that there are no $\Omega$-periodic $A$-modules. Finally, again by Corollary \ref{equivalencia periodicos}, we have that $\Omega^{\infty}(\md A) = \mathcal{P}(A)$. 

\end{itemize}

Finally consider $Q$ a quiver such that $Q^{\infty}$ has a final subheart equal to $C^h$ for $h \in \mathbb{Z}^+$.
Then there are two possible cases. 
\begin{itemize}

\item Every final subheart of $Q^{\infty}$ of type $C^h$ can be extended to a subquiver $\bar{Q}$ of $Q$ as below

$$
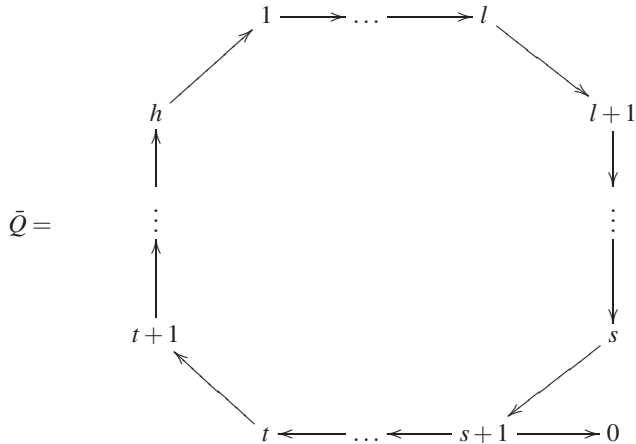
\xymatrix{ & & 1 \ar[r] & \ldots \ar[r] & l \ar[rd] &  \\
                  & h \ar[ru] & &  & & l+1 \ar[d] \\
                  \bar{Q} = & \vdots \ar[u] &  &   &  & \vdots \ar[d] \\
                  & t+1 \ar[u] &  &  & & s \ar[ld] \\ 
 	              & & t \ar[lu] & \ldots \ar[l] & s+1 \ar[l] \ar[r] & 0  \\}$$

\item There is a final subheart $Q^{\infty}$ of type $C^h$ that can not be extended as the previous case but can be extended to a subquiver $\bar{Q}$ of $Q$ as below

$$
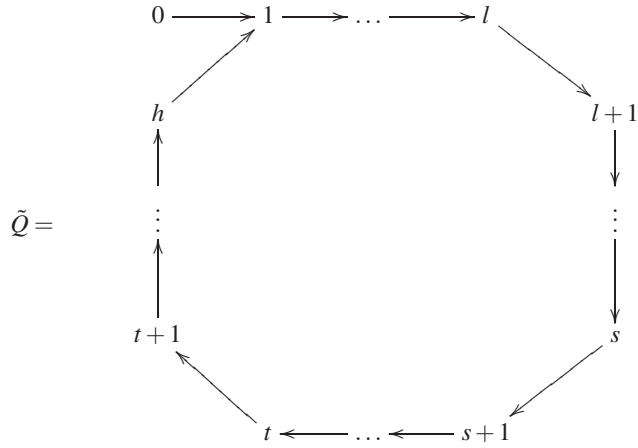
\xymatrix{& 0 \ar[r] & 1 \ar[r] & \ldots \ar[r] & l \ar[rd] &  \\
                  & h \ar[ru] & &  & & l+1 \ar[d] \\
                 \tilde{Q} =  & \vdots \ar[u] &  &   &  & \vdots \ar[d] \\
                  & t+1 \ar[u] &  &  & & s \ar[ld] \\ 
 	              & & t \ar[lu] & \ldots \ar[l] & s+1 \ar[l] &  \\}$$

\end{itemize}

In the second case, it is clear that there is a path $\rho: i \rightarrow j$ with $i, j \in \{1, \ldots , h \} $ such that the right ideal $\rho A$ is $A$-periodic but $\rho$ is not a perfect path in $Q$ because there are two paths $\gamma_1$ and $\gamma_2$ such that $\gamma_1 \rho \in J^m$ and $\gamma_2 \rho  \in J^m$. Therefore $\rho A$ is an $A$-periodic non Gorenstein projective module, and we deduce that $A$ is not Co-Gorenstein. 

In the first case, it is clear that there is a path $\rho: i \rightarrow j$ with $i, j \in \{1, \ldots , h \} $ such that the right ideal $\rho \frac{\Bbbk Q^{\infty}}{J^m}$ is $\frac{\Bbbk Q^{\infty}}{J^m}$-periodic but $\rho$ is not a perfect path in $Q$ because there are two paths $\lambda_1$ and $\lambda_2$ such that $\rho \lambda_1 \in J^m$ and $ \rho \lambda_2 \in J^m$. By Proposition \ref{periodico}, there is an $A$-module $M \in \mathcal{P}^{< \infty}(A)$ such that $\rho A \oplus M$ is $A$-periodic and it is not Gorenstein projective, thus we deduce that $A$ is not Co-Gorenstein.\end{proof}

\end{teo}

\begin{obs} Let $A$ be an Artin algebra that verifies

\begin{itemize}

\item $A$ is a truncated path algebra,

\item $A$ is a Co-Gorenstein algebra, and

\item $A$ is not a Gorenstein algebra,

\end{itemize}

then $A$ is CM-free.  

\end{obs}

In \cite{C} the authors prove that radical square zero algebras are either self-injective or CM-free. The example below shows that the behaviour of monomial algebras is more complicated. It shows that there are Co-Gorenstein monomial algebras that are not Gorenstein nor CM-free.

\begin{ej} Consider the algebra $A = \frac{\Bbbk Q}{I}$, where

$$Q = \xymatrix{ 1 \ar@/^/[r]^{\alpha} & 2 \ar@/^/[l]^{\beta} \ar@(ur,dr)^{\gamma} } \mbox{ and } I = \langle  \beta \alpha, \gamma^2  \rangle$$

Note that $\pd (I_1\oplus I_2) = \infty$, hence $A$ is not a Gorenstein algebra. On the other hand, the module $\gamma A $ verifies $\Omega (\gamma A  ) \cong \gamma A$, and it is the only indecomposable periodic non projective module. In fact $\gamma A$ is the only indecomposable non projective Gorenstein projective module by Theorem \ref{CSZ}. We conclude that $A$ is a Co-Gorenstein algebra and it is not CM-free . 

\end{ej}

\section{$\Omega^{\infty}$-infinite representation type and infinite $\phi$-dimension}

\subsection{Igusa-Todorov function}

The following lemma is a generalization of Proposition 3.9 of \cite{LMM}. The same proof works for it, just put $\mathcal{C}$ instead of $\md A$.

\begin{lema}\label{invariante}
Let $A$ be an Artin algebra. If $\mathcal{C}$ is a full subcategory of $\md A$ closed by direct sums and direct summands and $\Omega( \mathcal{C}) \subset \mathcal{C}$, then

$$\fidim(\mathcal{C}) = \min \{l:\bar{\Omega}{|}_{\Omega^l(K_0(\mathcal{C}))} \text{ is injective} \}. $$

\end{lema}

\begin{teo}\label{triangular}
Let $Q$ be a quiver and $C = \frac{\Bbbk Q}{I}$. Suppose there exist two disjoint full subquivers of $Q$, $\Gamma$ and $\bar{\Gamma}$ such that 

\begin{itemize}

\item  $Q_0 = \Gamma_0 \cup \bar{\Gamma}_0$, 

\item there is at least one arrow from $\Gamma_0$ to $\bar{\Gamma}_0$, and there are no arrows of $Q$ from $\bar{\Gamma}_0$ to $\Gamma_0$,

\item if $\alpha \in Q_1$ with $\st(\alpha) \in \Gamma_0$ and $\tg(\alpha) \in \bar{\Gamma}_0$ then $\beta \alpha = 0 $ for every arrow $\beta$.

\end{itemize}

Then $\fidim (C) \leq \fidim (A) + \fidim (B) + 1$ where $A = \frac{\Bbbk \Gamma}{I \cap \Bbbk \Gamma}$ and $B = \frac{\Bbbk \bar{\Gamma}}{I \cap \Bbbk \bar{\Gamma}}$.

\begin{proof}
We denote by $\partial A = \{v\in {\Gamma}_0: \exists \ \alpha \in Q_1 \text{ with } \st(\alpha) = v \text{ and } \tg(\alpha) \in \bar{\Gamma}_0\}$ and by  $ \mathcal{P}_{\partial A}$ the full subcategory of $\md A$ generated by the projective $A$-modules $P_v$ with $v \in \partial A$.  
First note that if $M \in \md C$ and

$$0 \rightarrow K \rightarrow P \rightarrow M \rightarrow 0$$

is an exact sequence in $\md C$ with $P$ a projective $C$-module, then $K \cong M_1 \oplus M_2 \oplus M_3 \oplus Q$ where $M_1 \in \Omega_A(\md A)$, $M_2 \in \mathcal{P}_{\partial A}$, $M_3\in \md B$ and $Q$ is a projective $C$-module.

\underline{Claim} $\phi (K) \leq \fidim (A) + \fidim (B)$

It is clear that
\begin{itemize}
\item  $\Omega_C(M_1) = \Omega_A (M_1) \oplus N_1 \oplus N_2 $ with $N_1 \in \mathcal{P}_{\partial A}$ and $N_2 \in \md B$,

\item $\Omega_C (M_2) \in \md B$, and

\item $\Omega_C (M_3) = \Omega_B (M_3)$.

\end{itemize}

If $\mathcal{D} = \Omega_A(\md A) \oplus \mathcal{P}_{\partial A} \oplus \md B$, then $K_0(\mathcal{D}) = K_1(A) \times K_0(\mathcal{P}_{\partial A}) \times K_0(B)$ and $\bar{\Omega}_C$ restricted to $K_0(\mathcal{D})$ is

$$ \left(
      \begin{array}{ccc}
        \bar{\Omega}_A & 0 & 0\\
        T_1 & 0 & 0\\
        T_2 & T_3 & \bar{\Omega}_B
      \end{array}
    \right).$$

Suppose that $\fidim (C) \geq k>0$, then, by Lemma \ref{invariante}, there exists $k' \geq k-1$ and $x = (x_1,x_2,x_3) \in K_0(\mathcal{D})$ such that $\bar{\Omega}_C^{k'}(x) = 0$ and $\bar{\Omega}_C^{k'-1}(x) \neq 0$. Since $\bar{\Omega}^i_A(x_1) =\Pi_{K_1(A)} (\bar{\Omega}^i_C(x))$ then $\bar{\Omega}^{k'}_A(x_1) = 0$.   We deduce that $\bar{\Omega}^l_A(x_1) = 0\ \forall l \geq \fidim (A) -1$ because $x_1 \in K_1(A)$, then $\bar{\Omega}_C^{\fidim (A) -1}(x) \in \{0\} \times K_0(\mathcal{P}_{\partial A}) \times K_0(B)$. Hence $k' \leq \fidim (A) + \fidim (B)$. Finally, the thesis follows by Proposition \ref{it1}.

\end{proof}

\end{teo}

\subsection{Examples}

 In \cite{M} it is showed that some Gorenstein algebras are of $\Omega^{\infty}$-infinite representation type (see example 4.1), hence  this algebra has finite $\phi$-dimension (see Theorem 4.7 of \cite{LM}). We now give an example of a non Gorenstein algebra of $\Omega^{\infty}$-infinite representation type, but with finite $\phi$-dimension.  

\begin{ej}\label{finito} Let $\Bbbk$ be an infinite field with characteristic different to $2$. Consider the algebra $A = \frac{\Bbbk Q}{•I}$, where $Q$ is the quiver below:

$$ Q = 
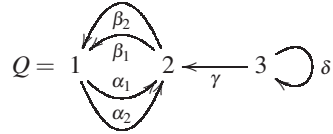
\xymatrix{ 1 \ar@/_{1pc}/[r]^{\alpha_1} \ar@/_{2pc}/[r]^{\alpha_2} & 2 \ar@/_{1pc}/[l]^{\beta_1} \ar@/_{2pc}/[l]^{\beta_2} & 3 \ar[l]^{\gamma} \ar@(ur,dr)^{\delta}}$$

and $$I = \langle \beta_1\alpha_1-2\beta_2\alpha_2,\ \alpha_1\beta_1 - \alpha_2\beta_2, \ \delta^2 , \ \delta\gamma , \ \gamma\beta_1, \gamma\beta_2 , \ \alpha_i \beta_j, \ \beta_i \alpha_j \mbox{ with } \ i,j \in \{1,2\} \mbox{ and } i \not = j  \rangle.$$

With the same computation as Example \ref{ej sizigia-infinita} we conclude that $A$ is of $\Omega^{\infty}$-infinite representation type.
It can be seen that the injective module $I(3)$ has infinite projective dimension, hence $A$ is not a Gorenstein algebra. 
Finally by Theorem \ref{triangular} we conclude that $\fidim(A) = 1$.

\end{ej}

In \cite{LM} is proved that $\Omega^n$-finite representation type algebras have finite $\phi$-dimension, therefore if we want to find an example of an algebra with infinite $\phi$-dimension, it must be of $\Omega^{\infty}$-infinite representation type.

The following example shows a radical cube zero algebra with infinite $\phi$-dimension. It is a counterexample of the conjecture stated in 
\cite{FLM}, however its finitistic dimension is finite (see \cite{Z2}). Since $\phi$DC is not valid, then Theorem \ref{FLM} in particular implies that the infinite $\phi$-dimensions are preserved by derived equivalences. The example also shows that the subcategory $\Omega^{\infty}(\md A)$ can be trivial although the algebra $A$ is of $\Omega^{\infty}$-infinite representation type.

\begin{ej}\label{infinito} Let $A = \frac{\Bbbk Q}{I}$ be an algebra where $Q$ is 

$$
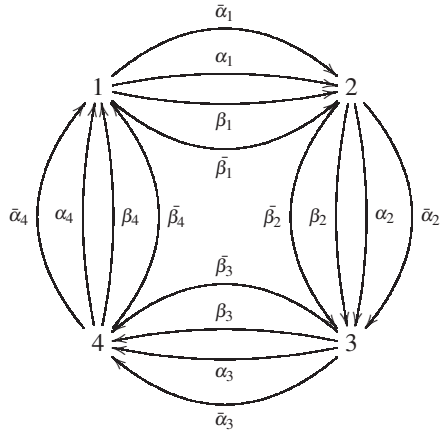
\xymatrix{ 1 \ar@/^8mm/[rrr]^{\bar{\alpha}_1} \ar@/^2mm/[rrr]^{\alpha_1} \ar@/_2mm/[rrr]_{\beta_1} \ar@/_8mm/[rrr]_{\bar{\beta_1}} & &  & 2 \ar@/^8mm/[ddd]^{\bar{\alpha}_2} \ar@/^2mm/[ddd]^{\alpha_2} \ar@/_2mm/[ddd]_{\beta_2} \ar@/_8mm/[ddd]_{\bar{\beta_2}} \\ & &  &\\& & & & \\ 4 \ar@/^8mm/[uuu]^{\bar{\alpha}_4} \ar@/^2mm/[uuu]^{\alpha_4} \ar@/_2mm/[uuu]_{\beta_4} \ar@/_8mm/[uuu]_{\bar{\beta_4}} &  & &  3 \ar@/^8mm/[lll]^{\bar{\alpha}_3} \ar@/^2mm/[lll]^{\alpha_3} \ar@/_2mm/[lll]_{\beta_3} \ar@/_8mm/[lll]_{\bar{\beta_3}}},$$

and $I = \langle \alpha_{i}\alpha_{i+1}-\bar{\alpha}_{i}\bar{\alpha}_{i+1},\ \beta_{i}\beta_{i+1}-\bar{\beta}_{i}\bar{\beta}_{i+1},\ \alpha_{i}\bar{\alpha}_{i+1},\ \bar{\alpha}_{i}\alpha_{i+1},\ \beta_{i}\bar{\beta}_{i+1},\ \bar{\beta}_{i}\beta_{i+1},\text{ for } i \in \mathbb{Z}_4, \ J^3 \rangle$

It is clear that $\pd(S_i) = \infty$ for every $i \in \mathbb{Z}_4$.
Let $M^{\alpha_1}_n$ and $M^{\beta_1}_n$, with $n \in \mathbb{Z}^+$, be the indecomposable $A$-modules defined by

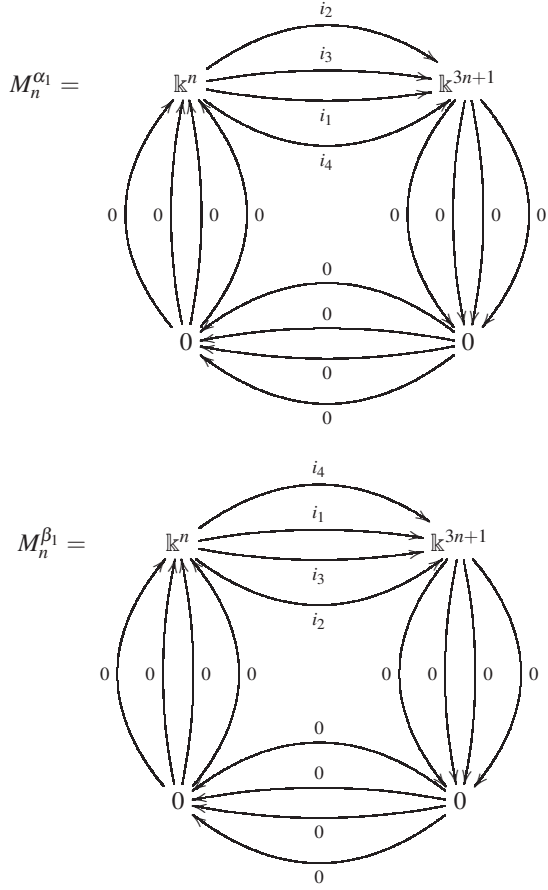
\[
M^{\alpha_1}_n =   \xymatrix{ & \Bbbk^n \ar@/^8mm/[rrr]^{i_2} \ar@/^2mm/[rrr]^{i_3} \ar@/_2mm/[rrr]_{i_1} \ar@/_8mm/[rrr]_{i_4} & &  & \Bbbk^{3n+1} \ar@/^8mm/[ddd]^{0} \ar@/^2mm/[ddd]^{0} \ar@/_2mm/[ddd]_{0} \ar@/_8mm/[ddd]_{0} \\ & & &  & \\ & & & & & \\ & 0 \ar@/^8mm/[uuu]^{0} \ar@/^2mm/[uuu]^{0} \ar@/_2mm/[uuu]_{0} \ar@/_8mm/[uuu]_{0} &  & &  0 \ar@/^8mm/[lll]^{0} \ar@/^2mm/[lll]^{0} \ar@/_2mm/[lll]_{0} \ar@/_8mm/[lll]_0},
\]

\[
 M^{\beta_1}_n = \xymatrix{ \Bbbk^n \ar@/^8mm/[rrr]^{i_4} \ar@/^2mm/[rrr]^{i_1} \ar@/_2mm/[rrr]_{i_3} \ar@/_8mm/[rrr]_{i_2} & &  & \Bbbk^{3n+1} \ar@/^8mm/[ddd]^{0} \ar@/^2mm/[ddd]^{0} \ar@/_2mm/[ddd]_{0} \ar@/_8mm/[ddd]_{0} \\ & &  &\\& & & & \\ 0 \ar@/^8mm/[uuu]^{0} \ar@/^2mm/[uuu]^{0} \ar@/_2mm/[uuu]_{0} \ar@/_8mm/[uuu]_{0} &  & &  0 \ar@/^8mm/[lll]^{0} \ar@/^2mm/[lll]^{0} \ar@/_2mm/[lll]_{0} \ar@/_8mm/[lll]_0},
\]

where the linear transformations $i_m: \Bbbk^n \rightarrow \Bbbk^{3n+1}$, with $m \in \{1,2,3,4\}$, verifies:

\begin{itemize}

\item  $i_1(e_j) = f_j\ \forall j \in \{1, \ldots n\}$,

\item  $i_2(e_j) =f_{n+j}\ \forall j \in \{1, \ldots n\}$,

\item  $i_3(e_j) =f_{n+j+1}\ \forall j \in \{1, \ldots n\}$,

\item  $i_4(e_j) = f_{2n+j + 1}\ \forall j \in \{1, \ldots n\}$.

\end{itemize}

where $\{e_1 \ldots e_n\}$ and $\{f_1,\ldots, f_{3n+1}\}$ are the canonical bases of $\Bbbk^n$ and $\Bbbk^{3n+1}$ respectively.

In an analogous way we define $M^{\alpha_2}_n, M^{\alpha_3}_n, M^{\alpha_4}_n$ and $M^{\beta_2}_n, M^{\beta_3}_n, M^{\beta_4}_n$. Then, the following assertions are valid for $i\in \mathbb{Z}_4$. 

\begin{itemize}

\item  $M^{\alpha_i}_n \ncong M^{\beta_i}_n $ if $n \geq 2$ and

\item  $M^{\alpha_i}_1 \cong M^{\beta_i}_1$.

\end{itemize}

We recall that the representation of $P_1$ is 

\[ P_1 =   \xymatrix{ & \Bbbk \ar@/^8mm/[rrr]^{u_2} \ar@/^2mm/[rrr]^{u_3} \ar@/_2mm/[rrr]_{u_1} \ar@/_8mm/[rrr]_{u_4} & &  & \Bbbk^{4} \ar@/^8mm/[ddd]^{v_2} \ar@/^2mm/[ddd]^{v_3} \ar@/_2mm/[ddd]_{v_1} \ar@/_8mm/[ddd]_{v_4} \\ & & &  & \\ & & & & & \\ & 0 \ar@/^8mm/[uuu]^{0} \ar@/^2mm/[uuu]^{0} \ar@/_2mm/[uuu]_{0} \ar@/_8mm/[uuu]_{0} &  & &  \Bbbk^{10} \ar@/^8mm/[lll]^{0} \ar@/^2mm/[lll]^{0} \ar@/_2mm/[lll]_{0} \ar@/_8mm/[lll]_0},
\]

where $u_{j}(\bar{e}_1) = \bar{f}_j$ and $v_1$, $v_2$, $v_3$, $v_4$ verifies:
\begin{center}

\begin{tabular}{|c|c|c|c|c|}

\hline 

 & \(v_4\) &  \(v_1\) & \(v_3\) & \(v_2\) \\ 

\hline 

 \(\bar{f}_4\) & \(\bar{g}_{3}\) & 0 & \(\bar{g}_1\) & \(\bar{g}_2\) \\ 

\hline 

 \(\bar{f}_1\) & \(0\) & \(\bar{g}_3\) & \(\bar{g}_4\) & \(\bar{g}_{5}\) \\ 

\hline 

 \(\bar{f}_3\) & \(\bar{g}_6\) & \(\bar{g}_7\)  & \(\bar{g}_8\) & 0 \\ 

\hline 

 \(\bar{f}_2\) & \(\bar{g}_9\) & \(\bar{g}_{10}\) & 0 &  \(\bar{g}_8\) \\ 

\hline 

\end{tabular} 
\end{center}

where $\{\bar{e}_1\}$, $\{\bar{f}_1, \bar{f}_2, \bar{f}_3, \bar{f}_4 \}$ and $\{\bar{g}_1, \bar{g}_2, \ldots, \bar{g}_9, \bar{g}_{10}\}$ are the canonical bases of $\Bbbk$, $\Bbbk^4$ and $\Bbbk^{10}$ respectively. $P_2$, $P_3$ and $P_4$ are analogous.

Since $\Tp (M^{\alpha_i}_n) = (S_i)^n$, the projective cover $M^{\alpha_i}_n$ is $(P_i)^n$.
Consider $i = 1$ (the other cases are analogous).

Let $F = (F_i)_{i \in 1, \ldots, 4}: (P_1)^n \rightarrow M^{\alpha_1}_n$ be the representation map associated a the projective cover that verifies $F_1((\bar{e}_1)_j) = e_j$ for $j=1, \ldots, m$. Notice that a representation for $\Omega(M^{\alpha_1}_n) = \Ker F$ is 

\[ \Omega(M^{\alpha_1}_n)  =   \xymatrix{ & 0 \ar@/^8mm/[rrr]^{0} \ar@/^2mm/[rrr]^{0} \ar@/_2mm/[rrr]_{0} \ar@/_8mm/[rrr]_{0} & &  & \Bbbk^{n-1} \ar@/^8mm/[ddd]^{w_2} \ar@/^2mm/[ddd]^{w_3} \ar@/_2mm/[ddd]_{w_1} \ar@/_8mm/[ddd]_{w_4} \\ & & &  & \\ & & & & & \\ & 0 \ar@/^8mm/[uuu]^{0} \ar@/^2mm/[uuu]^{0} \ar@/_2mm/[uuu]_{0} \ar@/_8mm/[uuu]_{0} &  & &  \Bbbk^{10n} \ar@/^8mm/[lll]^{0} \ar@/^2mm/[lll]^{0} \ar@/_2mm/[lll]_{0} \ar@/_8mm/[lll]_0},
\]
 
where the basis of $\Bbbk^{n-1}$ (at the vertex $2$) is $\{(-1)^i((f_3)_i - (f_2)_{i+1})\}_{i=1, \ldots n_1}$ and of $\Bbbk^{10n}$ (at the vertex $3$) is $\{(g_j)_i\}_{j = 1, \ldots 10, i=1, \ldots n }$.\\

The previous representation verifies

\begin{itemize}
\item The linear maps $w_1$, $w_2$, $w_3$ and $w_4$ are injective. 
\item The sum $w_1(\Bbbk^{n-1}) \oplus w_4(\Bbbk^{n-1})\oplus(w_2(\Bbbk^{n-1})+ w_3(\Bbbk^{n-1}))$
\item $w_2((-1)^i((f_3)_i-(f_2)_{i+1})) = (-1)^{i+1}(g_8)_{i+1} = w_3((-1)^{i+1}((f_3)_{i+1}-(f_2)_{i+2}))$
\end{itemize}

From the facts above we deduce that the syzygies of the modules $M^{\alpha_i}_n$ and $M^{\beta_i}_n$  for $n \geq 2$ and $i\in \mathbb{Z}_4$ are

\begin{itemize}

\item $\Omega(M^{\alpha_i}_n) = M^{\alpha_{i+1}}_{n-1} \oplus (S_{i+2})^{7n+2}$,

\item $\Omega(M^{\beta_i}_n) = M^{\beta_{i+1}}_{n-1} \oplus (S_{i+2})^{7n+2}$.

\end{itemize}

And for $n = 1$, we have $\Omega(M^{\alpha_i}_1) = \Omega(M^{\beta_i}_1) = (S_{i+2})^{10}$. 

By the previous facts we have $\phi(M^{\alpha_i}_n \oplus M^{\beta_i}_n) = n-1$ for any $n \geq 2$, then we conclude that $\fidim(A) = \infty$.

\end{ej}

Another counterexample to the $\phi$-dimension conjecture was given independently in \cite{HI}. Both examples are different since the algebra $A$ from Example \ref{infinito} is not a (nontrivial) tensor product of algebras (see the following remark). 

\begin{obs}
Let $A = \frac{\Bbbk Q}{I}$ be an algebra with $\vert Q_0 \vert = 4$. If $A = \frac{\Bbbk \tilde {Q}}{\tilde{I}} \otimes \frac{\Bbbk \bar{Q}}{\bar{I}}$ with $\tilde{Q}_1$ and $\bar{Q}_1$ non empty, then we have two different cases.

\begin{itemize}
\item $\vert \bar{Q}_0 \vert = 4$, $\vert \tilde{Q}_0\vert = 1$ and $\tilde{Q}$ has a loop.  

\item $\vert\tilde{Q}_0\vert = \vert \bar{Q}_0 \vert = 2$
\end{itemize}

In the first case $Q$ must have a loop. In the second one $Q$ has three different vertices $v, w, w'$ such that there are two arrows $\alpha, \beta \in Q_1$ such that $\st(\alpha) = \st(\beta) = v$, $\tg(\alpha) = w$ and $\tg(\beta) = w'$.

\end{obs}
 
\begin{obs}

Let $A$ be the algebra defined in Example \ref{infinito}, then every $A$-module $M \in  \Omega^1(\md A)$ verifies
$$\dim (\Tp(\Omega(M))) \geq \dim (\Tp(M))+1.$$ 

Consider an indecomposable $A$-module $M \in \Omega^1(\md A)$. We can suppose that $\underline{\dim}(M) = (\alpha, \beta, 0, 0)$, where $1 \leq \beta \leq 4\alpha-1$. Then $\underline{\dim}(\Omega(M)) = (0, 4\alpha-\beta, 10 \alpha, 0)$.

If $4\alpha-\beta \geq \alpha + 1$, we have $\dim (\Tp(\Omega(M))) \geq \dim (\Tp(M))+1$. Now suppose that $4\alpha-\beta \leq \alpha$, then $(S_3)^{6n}$ is a direct summand of $\Omega (M)$ therefore  $\dim (\Tp(\Omega(M))) \geq 6\dim (\Tp(M))$.

As a consequence of the above computation we conclude that $\Omega^{\infty}(\md A) = \mathcal{P}(A)$. Therefore $A$ is a CM-free Co-Gorenstein algebra.

\end{obs}

From $Corollary$ 3.5 of \cite{W} we conclude that $A$ of Example \ref{infinito} is $1$-Igusa Todorov.


\begin{thebibliography}{00}

\bibitem{B} A. Beligiannis. {\em The homological theory of contravariantly finite subcategories:
Auslander-Buchweitz contexts, Gorenstein categories and (co-)stabilization}. Comm. Algebra, {\bf 28} (10) pp. 4547-4596 (2000).

\bibitem{BMR} M. Barrios, G. Mata, G. Rama. {\em{Igusa-Todorov $\phi$ function for truncated path algebras}}, Algebras and Representation Theory {\bf 23} (3), pp. 1051-1063 (2020). 

\bibitem{BZR} E. Babson, B. Huisgen-Zimmermann, and R. Thomas, {\em Generic representation theory of
quivers with relations}, J. Algebra {\bf{322}} (6), pp. 1877-1918 (2009).

\bibitem{CSZ} X. Chen, D. Shen, G. Zhou. {\em{The Gorenstein-projective modules over a monomial algebra}}, Proceedings of the Royal Society of Edinburgh Section A: Mathematics {\bf 148} (6), pp. 1115-1134 (2018).

\bibitem{C} X. Chen, {\em{Algebras with radical square zero are either self-injective or CM-free}}, Proc. Am. Math. Soc. {\bf{140}} (1) pp. 93-98 (2012).

\bibitem{DH-ZL} A. Dugas, B. Huisgen-Zimmermann and J. Learned, {\em Truncated path algebras are homologically transparent. Part I}, Models, Modules and Abelian Groups (R. G\"{o}bel and B. Goldsmith, eds.), de Gruyter, Berlin, pp. 445-461 (2008).

\bibitem{FLM} S. Fernandes, M. Lanzilotta, O. Mendoza, {\em The Phi-dimension: A new homological measure}, Algebras and representation theory {\bf{17}} (5), pp. 463-476 (2014).

\bibitem{H} M. Hoshino, {\em Algebras of finite self-injective dimension}, Proc. Amer. Math. Soc.
{\bf 112} (3), pp. 619-622 (1991).

\bibitem{HI} E. Hanson, K. Igusa, {\em A Counterexample to the $\phi$-Dimension Conjecture}, arXiv:1911.00614 [math.RT] (2019).

\bibitem{HLM} F. Huard, M. Lanzilotta, O. Mendoza. {\em An approach to the finitistic dimension conjecture}, J. Algebra {\bf 319} (9), pp. 3916-3934 (2008).

\bibitem{IT} K. Igusa, G. Todorov, {\em On finitistic global dimension conjecture for artin algebras}, Representations of algebras and related topics, Fields Inst. Commun. {\bf 45}, American Mathematical Society, Providence, RI, pp. 201-204 (2005).

\bibitem{KM} S. Kvamme, R. Marczinzik, {\em{Co-Gorenstein algebras}}, Applied Categorical Structures {\bf{27}} (3), pp. 277-287 (2019).

\bibitem{M} G. Mata, {\em Igusa-Todorov function on path rings}, Bol. Soc. Paran. Mat. To appear.

\bibitem{LM} M. Lanzilotta, G. Mata, {\em Igusa-Todorov functions for Artin algebras}, Journal of pure and applied algebra, {\bf 222} (1), pp. 202-212 (2018). 

\bibitem{LMM} M. Lanzilotta, E. Marcos, G. Mata, {\em Igusa-Todorov functions for radical square zero algebras}, J. Algebra {\bf 487}, pp. 357-385 (2017). 

\bibitem{W} J. Wei, {\em Finitistic dimension and Igusa-Todorov algebras}, Advances in Mathematics {\bf 222} pp. 2215-2236 (2009).

\bibitem{Z} B. Zimmermann-Huisgen, {\em Predicting syzygies over monomial relations algebras}, Manuscripta Math. {\bf 70} (1), pp. 157-182 (1991).

\bibitem{Z2} B. Zimmermann-Huisgen, {\em Bounds on finitistic and global dimension for artinian rings with vanishing radical cube}, J. Algebra {\bf 161} pp. 47-68 (1993).

\bibitem{Zh} P. Zhang. {\em A brief introduction to Gorenstein projective modules}, Notes \url{https://www.math.uni-bielefeld.de/~sek/sem/abs/zhangpu4.pdf}

\end{thebibliography}
\end{document}